\documentclass[12pt]{article}
\usepackage{amssymb}
\newcommand{\sect}[1]{\setcounter{equation}{0}\section{#1}}

\newcommand{\be}{\begin{equation}}
\newcommand{\ee}{\end{equation}}
\newcommand{\bea}{\begin{eqnarray}}
\newcommand{\eea}{\end{eqnarray}}
\newcommand{\beano}{\begin{eqnarray*}}
\newcommand{\eeano}{\end{eqnarray*}}
\newcommand{\nonu}{\nonumber \\}

\newcommand{\hs}[1]{\hspace{#1 mm}}
\newcommand{\eps}{\epsilon}

\newcommand{\vph}{\varphi}
\newcommand{\sig}{\sigma}

\newcommand{\cc}{\mbox{$\cal{C}$}}

\newcommand{{\cg}}{\mbox{$\cal{G}$}}

\newcommand{\ci}{\mbox{${\cal I}$}}
\newcommand{\cj}{\mbox{${\cal J}$}}
\newcommand{\ck}{\mbox{$\cal{K}$}}
\newcommand{\cl}{\mbox{$\cal{L}$}}

\newcommand{\cu}{\mbox{${\cal U}$}}
\newcommand{\cw}{\mbox{$\cal{W}$}}

\newcommand{\wh}[1]{\widehat{#1}}

\newcommand{\mb}[1]{\hs{4}\mbox{#1}\hs{4}}
\newcommand{\half}{\frac{1}{2}}

\newcommand{\bn}{{\bar{n}}}
\newcommand{\sgn}{{\mbox{sg}}}
\newcommand{\qdet}{{\mbox{qdet}}}
\newcommand{\sdet}{{\mbox{sdet}}}

\newtheorem{prop}{Property}
\newtheorem{coro}{Corollary}
\newtheorem{defi}{Definition}
\newtheorem{theor}{Theorem}

\newcommand{\prf}{\underline{Proof:}\ }
\newcommand{\finprf}{\null \hfill {\rule{5pt}{5pt}}\\[2.1ex]\indent}
\newcommand{\ie}{{\it i.e.}\ }
\newcommand{\CC}{\mbox{${\mathbb C}$}}

\newcommand{\NN}{\mbox{${\mathbb N}$}}

\newcommand{\II}{\mbox{${\mathbb I}$}}

\newcommand{\FF}{\mbox{${\mathbb F}$}}

\newcommand{\NP}[1]{Nucl.\ Phys.\ {\bf #1}}
\newcommand{\PL}[1]{Phys.\ Lett.\ {\bf #1}}

\newcommand{\CMP}[1]{Commun.\ Math.\ Phys.\ {\bf #1}}

\begin{document}
\renewcommand{\thefootnote}{\fnsymbol{footnote}}
\newpage
\setcounter{page}{0}

\newcommand{\LAP}{LAPTH}
\def\logo{{\bf {\huge LAPTH}}}
\pagestyle{empty}

\centerline{\logo}

\vspace {.3cm}

\centerline{{\bf{\it\Large 
Laboratoire d'Annecy-le-Vieux de Physique Th\'eorique}}}

\centerline{\rule{12cm}{.42mm}}

\vspace{20mm}
\begin{center}
  {\LARGE  {\sffamily Twisted Yangians and folded \cw-algebras
    }}\\[1cm]

  
{\Large E. Ragoucy\footnote{ragoucy@lapp.in2p3.fr}}\\[.21cm] 

{\large Laboratoire de Physique Th{\'e}orique \LAP\footnote{UMR 5108 
du CNRS, associ{\'e}e {\`a} l'Universit{\'e} de Savoie.}\\[.242cm]
LAPP, BP 110\\
F-74941  Annecy-le-Vieux Cedex, France. }
\end{center}
\vfill\vfill

\begin{abstract}
 We show that the truncation of twisted Yangians are isomorphic to 
 finite \cw-algebras based on orthogonal or symplectic algebras. This
 isomorphism allows us to classify all the finite dimensional 
 irreducible representations of the quoted \cw-algebras. We also give 
 an R-matrix for these \cw-algebras, and determine their center.
\end{abstract}

\vfill
\rightline{\tt math.QA/0012182}
\rightline{\LAP-824/00}
\rightline{December 2000}

\newpage
\pagestyle{plain}
\setcounter{footnote}{0}
%
\pagestyle{plain}

\sect{Introduction}
Recently \cite{cmp,wrtt}, it has been remarked  that 
Yangians \cite{Drinfeld} and \cw-algebras \cite{zam} (both based on 
$gl(N)$ algebras) can be connected, 
although they belong to different fields of theoretical physics.

\cw-algebras have been introduced in the $2d$-conformal models as a 
tool for the study of these theories, such as the Toda field theories 
\cite{Oraf}. Then, these algebras and their 
finite-dimensional versions appeared to be relevant in several 
physical backgrounds \cite{Boscho,dBTj}. 
However, a full understanding of their algebraic 
structure (and of their geometrical interpretation) is 
lacking. 

Yangians were first considered and defined in
connection with some rational solutions of the quantum Yang-Baxter
equation. Later, their relevance in integrable models with non Abelian
symmetry was remarked \cite{4}.

The connection of some  of the finite \cw-algebras with Yangians (both based on 
$gl(N)$ algebras)
appears to shed some light on the algebraic structure of the former: it allows 
the construction of an $R$-matrix for \cw-algebras, the classification of their 
irreducible finite-dimensional representations and the determination of their center.

\null

In the present article, we continue with the study of this connection for 
the case of finite \cw-algebras based on orthogonal and symplectic algebras. 
Such \cw-algebras appear to be truncations of twisted Yangians \cite{Olsh,MNO}. 
Using this relation, we give an $R$-matrix formulation for the 
corresponding \cw-algebras. Contrarily to the $RTT$ formulation 
encountered in the case of $gl(N)$, it 
appears here to be of $ABCD$-algebras type \cite{maillet}, \ie an 
$RSR'S$ formulation. Since these 
algebras are not Hopf algebras, this 
confirms the remark that \cw-algebras seem to have {\it no} natural 
Hopf structure. The $RSR'S$ formulation allows to
classify the irreducible finite-dimensional 
representations of the \cw-algebras and to determine their center.

\null

The article is organized as follows: in section \ref{rappels}, we 
remind some notions on Yangians and finite \cw-algebras. In section 
\ref{sectYtw}, we present the construction that lead to twisted Yangians, as 
it was originally done in \cite{MNO}, and summarize some of their 
properties, such as the classification of finite-dimensional 
irreducible representation. In section \ref{sectWfold}, we apply the 
same type of procedure to finite \cw-algebras, to obtain what is known 
as folded \cw-algebras \cite{fold}. The comparison between these two 
objects is done in section \ref{results}, and is used to classify the 
finite-dimensional representations of \cw-algebras. We conclude in 
section \ref{concl}.
\pagebreak[3]

\sect{Yangians and \cw-algebras based on $gl(Np)$\label{rappels}}
We briefly present some known results on Yangians \cite{Drinfeld} and finite
\cw-algebras \cite{dBTj} that will be used in the following.
\subsection{The Yangian $Y(gl(N))\equiv Y(N)$\label{sectYN}}
One starts with the Yangian $Y(N)$ based on $gl(N)$. It is a Hopf algebra, 
the structure of which is contained in the relations 
(see for instance \cite{YN} and ref. therein for more details):
\bea
&&R_{12}(u-v)\, T_{1}(u)\, T_{2}(v) = T_{2}(v)\, T_{1}(u)\, R_{12}(u-v)
\label{rtt} \\
&& \Delta(T(u))=T(u)\otimes T(u)
\eea
We use here the usual convention on auxiliary spaces 
$T_{1}=T\otimes\II$  and $T_{2}=\II\otimes T$.
The generators $T^{ij}_{(n)}$ of the Yangian
are gathered in 
\be
  T(u)=\sum_{n=0}^\infty 
\sum_{i,j=1}^N u^{-n} T_{(n)}^{ij} E_{ij}=\sum_{n=0}^\infty u^{-n} T_{(n)}
=\sum_{i,j=1}^N T^{ij}(u) E_{ij}  \mbox{ with } T^{ij}_{(0)}=\delta^{ij}
\ ;\ T_{(0)}=\II
\ee
where $E_{ij}$ is the $N\times N$ matrix with 1 at position $(i,j)$, 
and $\II=\sum_{i=1}^NE_{ii}$.

The $R$-matrix is given by 
\be
R_{12}(x) = \II\otimes\II-\frac{1}{x}P_{12}\ ;\ P_{12}=\sum_{i,j=1}^N 
E_{ij}\otimes E_{ji}
\ee 
$P_{12}$ is the permutation of the two auxiliary spaces.

The relation (\ref{rtt}) is equivalent to the commutation relations:
\be
{[T_{(m)}^{ij},T_{(n)}^{kl}]} = \sum_{r=0}^{min(m,n)-1} \left(
T_{(r)}^{kj}T_{(m+n-r-1)}^{il}-T_{(r)}^{il}T_{(m+n-r-1)}^{kj}\right)
\ee
or keeping the auxiliary spaces:
\be
{[T_{(m)1},T_{(n)2}]} = \sum_{r=0}^{min(m,n)-1} \left(
P_{12}T_{(r)1}T_{(m+n-r-1)2}-T_{(r)2}T_{(m+n-r-1)1}P_{12}\right)
\ee

$(n)$ defines a natural gradation on the Yangian, and we will call level
the corresponding grade (\ie $T_{(n)}^{ij}$ and  $T_{(n)}$ are
said of level $n$).

The center $Z(N)$ of $Y(N)$ has been determined in \cite{MNO,Olsh}. 
It is generated by the 
quantum determinant:
\be
\qdet[T(u)]=\prod_{\sig\in S_N} (-1)^{\mbox{\scriptsize{sg}}(\sig)}
 T_{1\sig(1)}(u) T_{2\sig(2)}(u-1)\ldots T_{N\sig(N)}(u-N+1) 
 =1+\sum_{{n>0}}d_{n}u^{-n}
 \ee
 Denoting by $SY(N)$ the Hopf algebra $Y(N)/Z(N)$, we have:
 \be
 Y(N)\equiv Z(N)\otimes SY(N)\ \mbox{ and }\ SY(N)\equiv Y(sl(N))
 \ee
 where $Y(sl(N))$ is the Yangian based on $sl(N)$.
 
 The finite dimensional representations of $Y(N)$ has also been determined 
 in \cite{Dreval,tara}, see also \cite{Chari,Mol} for more 
details. 
\subsubsection{Classical version}
One can take the classical limit of the Yangian:
\be
R(x)=\II\otimes\II-\hbar\ r(x)\ ;\ T(u)=L(u)\ ;\ [\cdot,\cdot]=\hbar\{\cdot,\cdot\}
\ee
to get a Poisson Bracket
\be
\{L_{1}(u),L_{2}(v)\}= [r_{12}(u-v), L_{1}(u)L_{2}(v)]
\label{classY}
\ee
This procedure defines a classical (Poisson bracket) version of the Yangian, 
with an Abelian product
$L_{(m)}^{ij}\,L_{(n)}^{kl}=L_{(n)}^{kl}\,L_{(m)}^{ij}$.

\subsection{The $\cw(gl(Np),N.sl(p))\equiv\cw_p(N)$ algebra\label{sectWpN}}
It is defined as an Hamiltonian reduction of $gl(Np)$, considered as a Poisson algebra
(\ie with Poisson brackets). 
A basis of $gl(Np)$ (see \cite{wrtt} for more details) 
consists in generators $J_{jm}^{ab}$, with 
$-j\leq m\leq j$, $0\leq j\leq p$, and $a,b=1,...,N$, submitted to
\be
\{J_{jm}^{ab}, J^{cd}_{\ell,n}\} = 
    \sum_{r=\vert j-\ell\vert}^{ j+\ell}\sum_{s=-r}^r \left(
    \rule{0ex}{2.4ex}
    \delta^{bc} <j,m;\ell,n\vert r,s> J^{ad}_{r,s} - 
    \delta^{ad} <\ell,n;j,m\vert r,s> J^{cb}_{r,s} \right)
    \label{PBgln}
\ee
 $<j,m;\ell,n\vert r,s>$ are some Clebsch-Gordan like coefficients, 
 defined by
\[
    <j,m;\ell,n\vert r,s> = \frac{(-1)^s}{N\,\eta_{r}} tr\left( M^{ab}_{j,m}\cdot 
    M^{bc}_{\ell,n}\cdot M^{ca}_{r,-s}\right)\ \mbox{ with }
    \eta_{r}= (2r)! (r!)^2 {p+r\choose 2r+1}
\]
where $M_{jm}^{ab}$ are $(Np)\times(Np)$ matrices representing
$J_{jm}^{ab}$ in the fundamental of $gl(Np)$. They have being given in 
\cite{wrtt}, and we take the same normalizations:
\beano
M_{jm}^{ab} &=& E^{ab}\otimes
\left(\sum_{k=1}^{p-m} a_{j,m}^k E_{k,k+m}\right)
 \mbox{ ; } 
  a_{j,m}^k = \sum_{i=0}^{j-m} (-)^{i+j+m} {j-m\choose i} 
    a_{j,j}^{k-i}\mbox{ for } m\geq0 \\
M_{jm}^{ab} &=&  E^{ab}\otimes
\left(\sum_{k=1}^{p+m} a_{j,m}^k E_{k-m,k}\right) 
\mbox{ ; } a_{j,m}^k = \sum_{i=0}^{j-m} (-)^{i+j+m} {j-m\choose i} 
    a_{j,j}^{k-i-m}  \mbox{ for } m\leq0 \\
a_{j,j}^{k} &=& \frac{(k+j-1)! (p-k)!}{(k-1)! (p-k-j)!}
\eeano
where $E^{ab}$ are $N\times N$ matrices and $E_{k\ell}$ are 
$p\times p$ ones. Note that we have the properties \cite{wrtt}:
\bea
a_{j,-j}^{k} &=& a_{j,-j}^{0} = (-1)^j\ (2j)! \\
a_{j,m}^k &=& (-1)^{j+m} a_{j,m}^{N+1-k-m} \\
a_{j,m}^k &=& 0 \ \mbox{ for }\ |m|>j 
\eea

On $gl(Np)$ we impose a set of second class constraints
\be
\begin{array}{l}
J_{jm}^{ab}=0\ \mbox{ for }m<0,\ \forall j,a,b\mbox{ but }m=-1,\ j=1\\
J_{1,-1}^{ab}=\delta^{ab}\ \forall a,b
\end{array}
\label{const}
\ee 
which will be denoted by $\Phi=\{\vph_\alpha\}_{\alpha\in I}$ 
for convenience. 

The \cw-algebra is defined as the enveloping algebra of the 
generators $J_{jj}^{ab}\equiv W_{j}^{ab}$ equipped
with the Dirac brackets associated to the constraints (\ref{const}):
\be
\{X,Y\}_*\sim \{X,Y\}-\sum_{\alpha,\beta\in 
I}\{X,\vph_{\alpha}\}\cc^{\alpha\beta}\{\vph_{\beta},Y\}
\ \ \forall X,Y
\ee
where the matrix $\cc^{\alpha\beta}$ is the inverse of the matrix 
of constraints:
\be
\cc_{\alpha\beta}\sim\{\vph_\alpha,\vph_\beta\}
\ ;\ \cc^{\alpha\beta}\cc_{\beta\gamma}\sim\cc^\alpha_\beta
\ee
The symbol $\sim$ means that one has to apply the constraints on 
the right hand side of each expression 
{\em once the Poisson Brackets have been computed}.

\subsubsection{The $\cw_p(N)$ algebra in the Yangian basis\label{Ytronc}}
In \cite{wrtt,cmp}, it has been shown that the $\cw(gl(Np),N.sl(p))\equiv\cw_p(N)$
algebras are truncations of the Yangian $Y(N)$. They
are thus defined by the following relations\footnote{Note that, 
due to conflicting notations between \cw-algebras and Yangians, 
$W_{n}^{ij}$ corresponds to $T_{(n+1)}^{ij}$.}:
\bea
&&T(u)=\sum_{n=0}^p \sum_{i,j=1}^N\ u^{-n} T_{(n)}^{ij} \otimes E_{ij}
\ ;\ T_{(0)}^{ij}=\delta^{ij}\label{Tu-W}\\
&&T_1(u)=T(u)\otimes\II\ ;\ T_2(u)=\II\otimes T(u)
\ ;\ r(x)=\frac{1}{x}P_{12}\\
&& \{T_1(u),T_2(v)\} =[r_{12}(u-v),T_1(u)T_2(v)]
\eea
or equivalently
\be
\{T_{(m)}^{ij},T_{(n)}^{kl}\} = \sum_{r=0}^{min(m,n,p)-1} \left(
T_{(r)}^{kj}T_{(m+n-r)}^{il}-T_{(r)}^{il}T_{(m+n-r)}^{kj}\right)
\ee

There quantization becomes very simple in this basis. It reads:
\bea
&&R_{12}(u-v)\, T_{1}(u)\, T_{2}(v) = T_{2}(v)\, T_{1}(u)\, 
R_{12}(u-v)\\
&&R_{12}(x) = \II\otimes\II-\frac{1}{x}P_{12}
\eea
with the same definition (\ref{Tu-W}) of $T(u)$.

This connection with the Yangian $Y(N)$ allows to classify all the 
finite dimensional irreducible representations of $\cw_{p}(N)$ and
to determine its center \cite{wrtt}.

\sect{Twisted Yangians\label{sectYtw}}
\subsection{Presentation of $Y^\pm(N)$}
As the $gl(N)$ algebra, $Y(N)$ possesses automorphisms, and in the 
same way one 
can reconstruct the $so(n)$ and $sp(2n)$ algebras from the $gl(N)$ ones, one 
is tempted to reconstruct the Yangians based on orthogonal and 
symplectic algebras from $Y(N)$. However, the situation appears to be more 
delicate in the case of Yangians. 
In fact, although an algebraic structure similar to the one 
of orthogonal and symplectic Yangians can be achieved using 
automorphisms, the resulting algebra is not isomorphic to these 
Yangians. It is even not a Hopf algebra. The construction has been 
defined in \cite{Olsh,MNO}.

More precisely, the automorphisms we consider on $Y(N)$ are of the form
\bea
&&\tau(T(u))=T^t(-u) \mbox{ with } T^t(u)=\sum_{i,j=1}^N T^{ij}(u) E^t_{ij}
\mbox{ and } E^t_{ij}=\theta^{i}\theta^{j}E_{N+1-j,N+1-i}\nonu
&&\tau(T^{ij}_{(m)}) = (-1)^m\theta^{N+1-i}\theta^{N+1-j}\, T^{N+1-j,N+1-i}_{(m)}\ ;\ 
\theta^{i}=\pm1
\eea
Asking the automorphism $\tau$ to be of order 2 leads to the constraint
\be
\theta^{i}\theta^{N+1-i}=\theta_0\ \mbox{ with }\ \theta_{0}=\pm1 
\label{eq:tau2}
\ee
If $N=2n+1$, equation (\ref{eq:tau2}) for $i=n+1$ implies $\theta_0=1$.
We thus have the following conditions on $\theta_0$:
\be
\theta_0=1\ \mbox{ for }\ N=2n+1\ \mbox{ and } \theta_0=\pm1 \ \mbox{ for }\ 
N=2n
\ee
Each allowed values of the parameters $\theta^i$ determine an automorphism
$\tau$. However, only the values of $\theta_0$ are relevant for our purpose, as we
will see in the following.

Once $\tau$ is chosen, one defines the following generators in $Y(N)$
\be
S(u) = T(u)\tau(T(u))
\ee
From the relations (\ref{rtt}), one deduces (see \cite{MNO} for more 
details) 
\be
R_{12}(u-v)\, S_{1}(u)\, R'_{12}(u+v)\, S_{2}(v) = 
S_{2}(v)\, R'_{12}(u+v)\, S_{1}(u)\, R_{12}(u-v)
\label{rsrs}
\ee
where 
\bea
&&R'(x)=(\tau\otimes id)(R(x))=(id\otimes\tau)(R(x))=\II\otimes\II
-\frac{1}{x}Q_{12}\nonu
&&\mbox{ with }Q_{12}=\sum_{i,j=1}^N \theta^{i}\theta^{j}
E_{ij}\otimes E_{N+1-i,N+1-j}=(\tau\otimes id)(P_{12})
\eea
Note that this construction is a particular case of the 
$abcd$-algebras introduced in \cite{maillet}, and was indeed mentioned there 
as an example. The full structure of twisted Yangians has been 
studied in \cite{MNO}.

The above relation defines a subalgebra $Y^\tau(N)$ of the Yangian 
$Y(N)$. Looking at the generators of $Y^\tau(N)$
\be
S(u)=\sum_{i,j=1}^N u^{-n} S_{(n)}^{ij} E_{ij}=\sum_{n=0}^\infty u^{-n} S_{(n)}
=\sum_{i,j=1}^N S^{ij}(u) E_{ij} \mbox{ with }S_{(0)}^{ij}=\delta^{ij}
\ee
one can show that they obey the following relations:
\bea
{[S_{(m)}^{ij},S_{(n)}^{kl}]} &=& \sum_{r=0}^{min(m,n)-1} 
\left[\rule{0ex}{1.42em}
S_{(r)}^{kj}S_{(m+n-r-1)}^{il}-S_{(m+n-r-1)}^{kj}S_{(r)}^{il}
+\right.\label{com-S}\\
&&\hspace{1.2em}
\left. +(-)^{n+r}\theta_0\left(
\theta^{k}\theta^{j}S_{(r)}^{i,N+1-k}S_{(m+n-r-1)}^{N+1-j,l}
-\theta^{i}\theta^{l}S_{(m+n-r-1)}^{k,N+1-i}S_{(r)}^{N+1-l,j}\right)
\rule{0ex}{1.42em}\right]+\nonu
&&+\sum_{s=0}^{[\frac{m-2}{2}]}\theta_0\theta^{j}\theta^{i}\left(
S_{(m-2s-2)}^{k,N+1-i}S_{(n+2s)}^{N+1-j,l}
-S_{(n+2s)}^{k,N+1-i}S_{(m-2s-2)}^{N+1-j,l}
\rule{0ex}{1.2em}\right)\nonumber
\eea
which is equivalent to (\ref{rsrs}), as well as
\bea
{[S^{ij}(u),S^{kl}(v)]} &=& \frac{1}{u-v}\ 
\left(\rule{0ex}{1.42em}
S^{kj}(u)S^{il}(v)-S^{kj}(v)S^{il}(u)\right)
+\\
&& -\frac{\theta_{0}}{u+v}\left(
\theta^{k}\theta^{j}\ S^{i,N+1-k}(u)S^{N+1-j,l}(v)
-\theta^{i}\theta^{l}\ 
S^{k,N+1-i}(v)S^{N+1-l,j}(u)\rule{0ex}{1.42em}\right)+
\nonu
&&+\frac{\theta_{0}\theta^{j}\theta^{i}}{u^2-v^2}\,
\left(S^{k,N+1-i}(u)S^{N+1-j,l}(v)
-S^{k,N+1-i}(v)S^{N+1-j,l}(u)
\rule{0ex}{1.42em}\right)\nonumber
\eea
and also
\beano
{[S_{1}(u),S_{2}(v)]} &=& 
\frac{1}{u-v}\left(P_{12}S_{1}(u)S_{2}(v)-S_{2}(v)S_{1}(u)P_{12}\rule{0ex}{1.2em}\right)+\\
&&-\frac{1}{u+v}\left(S_{1}(u)Q_{12}S_{2}(v)-S_{2}(v)Q_{12}S_{1}(u)\rule{0ex}{1.2em}\right)+\\
&&+\frac{1}{u^2-v^2}\left(P_{12}S_{1}(u)Q_{12}S_{2}(v)-S_{2}(v)Q_{12}S_{1}(u)P_{12}\rule{0ex}{1.2em}\right)
\eeano

\begin{theor}\label{theo-eps}
    All the $\theta^i$ dependence, but $\theta_0$, can be removed. We thus have only
one (two) different twisted Yangian(s), corresponding to the (two) value(s)
 $\theta_0=1$ ($\theta_0=\pm1$) for $N=2n+1$ (for $N=2n$).
\end{theor}
\prf
We prove this by exhibiting a basis in which the $\theta^i$ dependence, 
but $\theta_0$,
has disappeared. In fact, the twisted Yangians are generated by
\be\begin{array}{ll}
\displaystyle
Y^\pm(2n) & S^{ij}(u)\ ;\ S^{i,N+1-j}(u)\ ;\ S^{N+1-i,j}(u)\ \
i,j=1,\ldots,n \\
Y^+(2n+1) & \left\{\begin{array}{l}
S^{ij}(u)\ ;\ S^{i,N+1-j}(u)\ ;\ S^{N+1-i,j}(u)\\[1.2ex]
S^{n+1,n+1}(u)\ ;\ S^{n+1,i}(u)\ ;\ 
S^{i,n+1}(u)\ ;\ 
i,j=1,\ldots,n \\
\end{array}\right.\end{array}
\ee
One then defines for $Y^\pm(2n)$:
\be
{J}^{ij}(u)=\theta^{i}\theta^{j}S^{ij}(u)\ ;\  
{K}^{ij}(u)= \theta^{i}S^{i,N+1-j}(u)\ ;\  
{\bar K}^{ij}(u)= \theta^{j}S^{N+1-i,j}(u)
\ee
which obey to $\theta$-free commutation relations (except the $\theta_0$
dependence). For instance:
\bea
{[{J}^{ij}(u),{J}^{kl}(v)]} &=& \frac{1}{u-v}\ 
\left(\rule{0ex}{1.42em}
{J}^{kj}(u){J}^{il}(v)-{J}^{kj}(v) {J}^{il}(u)\right)
+\\
&& -\frac{\theta_{0}}{u+v}\left(
{K}^{i,k}(u) {\bar K}^{j,l}(v)
-{K}^{k,i}(v){\bar K}^{l,j}(u)\rule{0ex}{1.42em}\right)+
\nonu
&&+\frac{\theta_{0}}{u^2-v^2}\,
\left({K}^{k,i}(u) {\bar K}^{j,l}(v)
- {K}^{k,i}(v) {\bar K}^{j,l}(u)
\rule{0ex}{1.42em}\right)\\
\eea
For $Y^+(2n+1)$, the redefinition is the same as before, plus for 
the remaining generators:
\be
{J}_{0}(u)=\theta^{n+1} S^{n+1,n+1}(u)\ ;\  {L}^{i}(u)=\theta^{i} S^{n+1,i}(u)\ ;\ 
 {\bar L}^{i}(u)=\theta^{n+1}\theta^{i} S^{i,n+1}(u)
\ee
The commutation relations are then free from any $\theta$'s.
\finprf

\begin{defi}
The twisted Yangians $Y^\pm(N)$ correspond to the following choices 
for $\tau$:
\be
\begin{array}{lll}
\mbox{ For }Y^+(N) \ : &\theta^{i}=1,\, \forall i & \ie \theta_0=1\\
\mbox{ For }Y^-(2n) \ : &\theta^{i}=\sgn(\frac{N+1}{2}-i),\, \forall i & \ie \theta_0=-1
\end{array} 
\ee
With these choices, the Lie subalgebra for $Y^\pm(N)$ is 
\be
 so(N)\, \subset\, Y^+(N)\ \mbox{ and }\
sp(2n) \,\subset\, Y^-(2n)
\ee
\end{defi}

Denoting by $Y^\pm_k(N)$ the subset of $Y^\pm(N)$ formed by the generators 
of level $k$, we have
\be
\mbox{dim}Y^\pm_{2k-1}(N)=\frac{N(N\mp1)}{2}\ ;\ 
\mbox{dim}Y^\pm_{2k}(N)=\frac{N(N\pm1)}{2},\ k=1,2,\ldots
\ee

In \cite{MNO}, it has been proven: 
\begin{prop}\label{Qinv-tau}
  The $Y(N)$-subalgebra $Y^\tau(N)$ generated by $S(u) = 
  T(u)\tau(T(u))$
is isomorphic to the algebra defined by the two following relations:
\bea
&&R_{12}(u-v)\, S_{1}(u)\, R'_{12}(u+v)\, S_{2}(v) = 
S_{2}(v)\, R'_{12}(u+v)\, S_{1}(u)\, R_{12}(u-v)
\label{com}\\
&&\tau[S(u)]=S(u)+\theta_{0}\,\frac{S(u)-S(-u)}{2u}\label{tauS}
\eea
\end{prop}

Let us remark that one can perform the change:
\be
S'(u)=S(u)+\theta_{0} \frac{S(u)-S(-u)}{4u}
=\frac{1}{2}\left(\rule{0ex}{2.3ex}S(u)+\tau\left[S(u)\right]\right) 
\ee
which obeys $\tau[S'(u)]=S'(u)$. This proves that twisted Yangians are 
a subalgebra of Ker$(id-\tau)$, the subalgebra generated by 
elements of $Y(N)$ which are invariant under $\tau$. However, the 
commutation relations satisfied by $S'(u)$ are more complicated than 
the ones obeyed by $S(u)$. We will come back on this point in the 
classical case (see below).

\subsection{Center of $Y^\pm(N)$}
The center of the twisted Yangian has been studied in \cite{Olsh,MNO}.
It is a true subalgebra of the center of $Y(N)$, and is generated by 
the so-called Sklyanin determinant 
sdet. The exact expression of the Sklyanin determinant in term 
of $S(u)$ can be found in \cite{Mol3}. It is however a rather complicated 
expression.
A more easy-to-handle expression, which refers to the underlying $Y(N)$, 
can be found in \cite{Olsh,MNO}:
\be
\begin{array}{l}
\displaystyle
\sdet[S(u)] =\gamma_N(u)\, \qdet[T(u)]\, \qdet[T(N-1-u)]
\ \mbox{ with :}\\[2.1ex]
\displaystyle\gamma_N(u)=1 \ \mbox{ for }\ Y^+(N) \mb{ and }
\gamma_N(u)=\frac{2u+1}{2u+1-N} \ \mbox{ for }\ Y^-(N)
\end{array}
\ee
One can show \cite{Mol3,MNO} that the Sklyanin determinant provides a 
basis $c_{2},c_{4},\ldots$ for the center of $Y^\pm(N)$ through the formulae:
\be
\begin{array}{l}
\displaystyle
\sdet[T(u-\frac{N+1}{2})]=1+\sum_{n>0}c_{2n}u^{-2n}\ \mbox{ for }\ 
Y^+(N)\\[2.1ex]
\displaystyle
\sdet[T(u+n-\frac{1}{2})]=(1+n\, u^{-1})\left(1+\sum_{n>0}c_{2n}u^{-2n}
\right)\ \mbox{ for }\ Y^-(2n)
\end{array}
\label{center}
\ee

If one denotes by $Z^\pm(N)$ the center of $Y^\pm(N)$, one has 
$Y^\pm(N)\equiv Z^\pm(N)\otimes SY^\pm(N)$, where the special twisted 
Yangian $SY^\pm(N)$ is defined by $SY^\pm(N)=SY(N)\cap Y^\pm(N)$.

\subsection{Classical case\label{classical}}
As for the Yangian, one can take a classical limit of the twisted Yangian
\be
\{S_{1}(u),S_{2}(v)\} = [r_{12}(u-v), S_{1}(u)S_{2}(v)] +
S_{2}(v)r'_{12}(u+v)S_{1}(u)-S_{1}(u)r'_{12}(u+v)S_{2}(v)
\label{classYtw}
\ee
where 
$r'_{12}(x)=(id\otimes\tau)(r(x))=(\tau\otimes id)(r(x))=\frac{1}{x}Q_{12}$, or 
more explicitly
\beano
\{S_{1}(u),S_{2}(v)\} &=&\frac{1}{u-v}\left(\rule{0ex}{2.3ex}
P_{12}\, S_{1}(u)\,S_{2}(v)-
S_{2}(v)\,S_{1}(u)\,P_{12}\right)+\\
&&-\frac{1}{u+v}\left(\rule{0ex}{2.3ex}
S_{1}(u)\,Q_{12}\,S_{2}(v)-S_{2}(v)\,Q_{12}\,S_{1}(u)\right)
\eeano
It leads to the following Poisson brackets for the generators:
\bea
\{S_{(m)}^{ij},S_{(n)}^{kl}\} &=& \sum_{r=0}^{min(m,n)-1} 
\left[\rule{0ex}{1.42em}
S_{(r)}^{kj}S_{(m+n-r-1)}^{il}-S_{(m+n-r-1)}^{kj}S_{(r)}^{il}+
\right.\label{PB-S}\\
&&+(-)^{n+r}\theta_{0}\left.\left(
\theta^{k}\theta^{j}S_{(r)}^{i,N+1-k}S_{(m+n-r-1)}^{N+1-j,l}
-\theta^{i}\theta^{l}S_{(m+n-r-1)}^{k,N+1-i}S_{(r)}^{N+1-l,j}
\right)\rule{0ex}{1.42em}\right]\nonumber
\eea

Note that (\ref{classYtw}) can be defined as the classical limit of the 
twisted Yangian as well as the twisted subalgebra of the classical 
Yangian, with still $S(u)=L(u)\tau(L(u))$.

Of course, as in the quantum case, it is only the $\theta_0$ dependence which 
is relevant, and we will deal with classical twisted Yangian $Y^\pm(N)$ only.

\begin{prop}\label{inv-tau}
  The classical $Y^\tau(N)$ algebra, generated by $S(u) = 
  T(u)\tau(T(u))$,
is isomorphic to the $Y(N)$-subalgebra defined by the two following relations:
\beano
&&\{S_{1}(u),S_{2}(v)\} = [r_{12}(u-v), S_{1}(u)S_{2}(v)] +
S_{2}(v)r'_{12}(u+v)S_{1}(u)-S_{1}(u)r'_{12}(u+v)S_{2}(v)
\\
&&\tau(S(u))=S(u)
\eeano
The classical twisted Yangians $Y^\tau(N)$ are thus isomorphic to subalgebras of the 
$\tau$-invariant subalgebra Ker$(id-\tau)$ in $Y(N)$.
\end{prop}
\prf
We start with the classical version of the property \ref{Qinv-tau} (proven in 
\cite{MNO}): $Y^\pm(N)$ is completely defined by the relation 
(\ref{classYtw}), together with
\be
\tau[S(u)]=S(u)+\theta_{0}\,\frac{S(u)-S(-u)}{2u}
\ee
Now, defining 
\be
S'(u)=S(u)\pm 
\frac{S(u)-S(-u)}{4u}=\frac{1}{2}\left(\rule{0ex}{2.3ex}S(u)+
\tau\left(S(u)\right)\right) 
\ee
one computes that $S'(u)$ still obey the quadratic relation (\ref{classYtw}) 
with now as symmetry relation $\tau(S'(u))=S'(u)$. Indeed, starting 
from the quadratic relation on $S(u)$, and denoting $\bar{S}=\tau(S)$, one deduces:
\beano
\{\bar{S}_{1}(u),S_{2}(v)\} &=& [r_{12}(u-v), \bar{S}_{1}(u)S_{2}(v)] +
S_{2}(v)r'_{12}(u+v)\bar{S}_{1}(u)-\bar{S}_{1}(u)r'_{12}(u+v)S_{2}(v) \\
\{S_{1}(u),\bar{S}_{2}(v)\} &=& [r_{12}(u-v), S_{1}(u)\bar{S}_{2}(v)] +
\bar{S}_{2}(v)r'_{12}(u+v)S_{1}(u)-S_{1}(u)r'_{12}(u+v)\bar{S}_{2}(v)
\\
\{\bar{S}_{1}(u),\bar{S}_{2}(v)\} &=& [r_{12}(u-v), \bar{S}_{1}(u)\bar{S}_{2}(v)] +
\bar{S}_{2}(v)r'_{12}(u+v)\bar{S}_{1}(u)-\bar{S}_{1}(u)r'_{12}(u+v)\bar{S}_{2}(v)
\eeano
where one has heavily used the commutativity of the product at the classical 
level. From these formulae, it is simple to deduce that 
$S'=\half(S+\bar{S})$ also obeys (\ref{classYtw}).

Finally, the calculation
\be
S'(u)-S'(-u)=S(u)-S(-u)
\ee
shows that the change from $S(u)$ to $S'(u)$ is invertible.
\finprf

In the basis presented in theorem \ref{theo-eps}, the condition 
$\tau[S(u)]={S}(u)$ reduces (at classical level)
to ${K}^{ij}(u)=\theta_0 {K}^{ji}(-u)$, 
${\bar K}^{ij}=\theta_0{\bar K}^{ji}(-u)$, and ${J}_0(u)={J}_0(-u)$.
It is thus also independent from $\theta$ (except $\theta_0$).

In the following, when dealing with classical twisted Yangian, we will 
choose as generating system the one given in property \ref{inv-tau}.

\subsection{Representations of twisted Yangian}
Irreducible representations of twisted Yangians have been studied in 
\cite{Mol}.  We recall here some of the obtained results\footnote{Note 
the difference of convention: $i,j=1,\ldots,N$ (used here) with respect to  
 $-n\leq i,j\leq n$ used in \cite{Mol}. The correspondence is given by 
$i\ \rightarrow\ \bar{n}-i$ for $i<0$, and $i\ \rightarrow\ n+1-i$
for $i\geq0$. 
Hence, the apparition of $i=0$
in \cite{Mol} is associated to $n\neq\bar{n}$ (which both occur only when $N=2n+1$).
The transformation $i\ \rightarrow\ -i$ used in \cite{Mol} is 
 translated into $i\ \rightarrow\ N+1-i$ used in the present article.}.
In the following we use the notation
\be
n=\left[\frac{N}{2}\right]\ \mbox{ and }\ \bar{n}=\left[\frac{N+1}{2}\right]
\ee
Remark that $n=\bar{n}$ iff $N$ is even, and that $N=n+\bar{n}$ in all 
cases.
\subsubsection{Classification \label{Ytwrep}}
\begin{defi}
 A representation $V$ of $Y^\pm(N)$ is called lowest weight if there 
 exists a vector $\xi\in V$ such that
\beano
&& S^{ij}(u)\xi=0\ \ \forall\ 1\leq j<i\leq N\\
 && S^{i,i}(u)\xi=\mu^i(u)\xi\ \ \forall\ 1\leq i\leq 
\bar{n}
\eeano
$\xi$ is the lowest weight vector of $V$ and 
$\mu(u)=(\mu^{1}(u),\dots,\mu^{\bar n}(u))$  its lowest weight.
\end{defi}
\begin{prop}
 Any finite dimensional irreducible representation  of $Y^\pm(N)$ 
 is lowest weight. The lowest weight vector of such a
 representation is unique (up to multiplication by a scalar).
    
 One can chose the basis of  $Y^\pm(N)$  in such a way that we have
\bea
    \mu^{i}(u) &=& 
    \prod_{k=1}^{d_{i}}\left(1-\lambda^{i}_{(k)}u^{-1}\right)\ 
     \mbox{ for } Y^-(2n)\\
    (1+\half u^{-1})\,\mu^{i}(u) &=& \prod_{k=1}^{2d_{i}+1}
    \left(1-\lambda^{i}_{(k)}u^{-1}\right)\ 
    \mbox{ for }Y^+(2n)\\
    \mu^{i}(u) &=& \prod_{k=1}^{d_{i}}\left(1-\lambda^{i}_{(k)}u^{-1}\right)\ 
    \mbox{ for } Y^+(2n+1)
\eea
    where $d_{i}\in\NN$ and $\lambda^{i}_{(k)}\in\CC$.
\end{prop}
\begin{prop}\label{irrepsYtw}
    There is a one-to-one correspondence between finite-dimensional 
    representations of $Y^\pm(N)$ and the families 
    $\{P_{1}(u),\dots,P_{\bar n}(u),\rho(u),\eps\}$ where 
    $P_{i}(u)$ ($\forall\ i$) are 
    monic polynomials in $u$ with $P_{n}(u)=P_{n}(1-u)$, $\rho(u)$ is 
    a formal series in $u^{-1}$ which encodes the values of the Casimir 
    operators of 
    $Y^\pm(N)$ and the parameter $\eps$ can take the following values:
    \begin{itemize}
    \item $\eps=1$ for $Y^-(2n)$ (\ie $\eps$ is not relevant in this case).
    \item $\eps\in\CC$ for $Y^+(2)$, with $P_{1}(-\eps)\neq0$
    \item $\eps=$1,2,3 or 4 for $Y^{+}(2n)$, $n>1$, with the restriction that the 
     values $\eps=1$ and $\eps=3$ have to be  identified when
     $P_{n}(\frac{1}{2})\neq0$.
    \item $\eps=1$ or $2$ for $Y^{+}(2n+1)$.
\end{itemize}
\end{prop}
By monic polynomials we mean a polynomial of the form 
$P(u)=\prod_{k=1}^m(u-\gamma_{k})$ for some complexes $\gamma_{k}$ and 
some integer $m$. The monic polynomials are related to the lowest 
weight of the representation by the relations 
\be
\frac{\mu^{i+1}(u)}{\mu^{i}(u)}=\frac{P_{i}(u+1)}{P_{i}(u)},\ 
\forall\ i=1,\ldots,n-2
\ee
Together with two supplementary relations which depend of the studied 
twisted Yangian. 

In the case of $Y^-(2n)$, the conditions take the form
\be
\frac{\mu^{n}(-u)}{\mu^{n}(u)}=\frac{P_{n}(u+1)}{P_{n}(u)} 
\ \mbox{ and }\  
\frac{\mu^{n}(u)}{\mu^{n-1}(u)}=\frac{P_{n}(u+1)}{P_{n}(u)}
\ee
The case $Y^+(2)$ is special, the Lie subalgebra being the Abelian 
$o(2)$. The lowest weight is reconstructed using:
\be
\frac{\mu(u)}{\mu(-u)}=\frac{(2u+1)(u+\eps)}{(2u-1)(u-\eps)}\frac{P(u+1)}{P(u)} 
\ee
When considering $Y^+(2n)$, $n>1$, one has to choose one of 
the four following possibilities, the choices being labeled by 
the parameter $\eps$:
\beano
\eps=1&:&\frac{\mu^{n}(-u)}{\mu^{n}(u)}=\frac{P_{n}(u+1)}{P_{n}(u)}
\ \mbox{ and }\ 
\frac{\mu^{n}(u)}{\mu^{n-1}(u)}=\frac{P_{n}(u+1)}{P_{n}(u)}\\
\eps=2&:&\frac{2u-1}{2u+1}\ 
\frac{\mu^{n}(-u)}{\mu^{n}(u)}=\frac{P_{n}(u+1)}{P_{n}(u)}
\ \mbox{ and }\ 
\frac{\mu^{n}(u)}{\mu^{n-1}(u)}=\frac{P_{n}(u+1)}{P_{n}(u)}\\
\eps=3&:&\frac{\mu_{\#}^{n}(-u)}{\mu_{\#}^{n}(u)}=\frac{P_{n}(u+1)}{P_{n}(u)}
\ \mbox{ and }\ 
\frac{\mu_{\#}^{n}(u)}{\mu^{n-1}(u)}=\frac{P_{n}(u+1)}{P_{n}(u)}\\
\eps=4&:&\frac{2u-1}{2u+1}\ \frac{\mu_{\#}^{n}(-u)}{\mu_{\#}^{n}(u)}=
\frac{P_{n}(u+1)}{P_{n}(u)}\ \mbox{ and }\ 
\frac{\mu_{\#}^{n}(u)}{\mu^{n-1}(u)}=\frac{P_{n}(u+1)}{P_{n}(u)}
\eeano
where $\mu_{\#}^{n}(u)$ is deduced from $\mu^{n}(u)$ by the equations
\bea
(1+\half u^{-1})\ \mu^{n}(u) &=& \prod_{k=1}^{2d+1}\left(1-\lambda^{(k)}u^{-1}
\right)\\
(1+\half u^{-1})\ \mu^{n}_{\#}(u) &=& \left(1+(\lambda^{(2d+1)}+1)u^{-1}\right)
\prod_{k=1}^{2d}\left(1-\lambda^{(k)}u^{-1}\right)
\eea
Finally, for $Y^+(2n+1)$, one of the two following choices has to be 
selected, corresponding to the two possible values of $\eps$:
\beano
\eps=1&:&\frac{\mu^{\bar n}(u)}{\mu^{n}(u)}=\frac{P_{\bar n}(u+1)}{P_{\bar n}(u)}
\ \mbox{ and }\ 
\frac{\mu^{n}(u)}{\mu^{n-1}(u)}=\frac{P_{n}(u+1)}{P_{n}(u)}
\\
\eps=2&:&\frac{2u}{2u+1}\frac{\mu^{\bar n}(u)}{\mu^{n}(u)}=
\frac{P_{\bar n}(u+1)}{P_{\bar n}(u)}
\ \mbox{ and }\  
\frac{\mu^{n}(u)}{\mu^{n-1}(u)}=\frac{P_{n}(u+1)}{P_{n}(u)}
\eeano

\subsubsection{Construction from $Y(N)$ representations}
All the irreducible finite-dimensional representations of the twisted 
Yangians can be constructed starting from $Y(N)$ representations. 
Here, we summarize the construction and
 refer to \cite{Mol} for the complete presentation. 
The basic idea is to try to build a $Y^\pm(N)$ representation as the 
restriction of an $Y(N)$ one (and eventually doing the subquotient 
to get an irreducible representation). This provides irreducible 
representations of $Y^\pm(N)$, but to get a complete classification, 
one has (sometimes) to add 
 an $o(N)$ representation, understood as an evaluation 
representation of $Y^+(N)$.

Let $\mu(u)$ be a lowest weight entering in the classification of 
section \ref{Ytwrep}, and $V(\mu)$ the corresponding representation. 
Here, we give (without calculations, see \cite{Mol} for more details)
the lowest weight $\lambda(u)$ and the
irreducible representation $L(\lambda)$
of $Y(N)$ which lead to $V(\mu)$ in term of the polynomials $P_{i}(u)$ 
of property \ref{irrepsYtw}. 

{\bf We start with the case of $\mathbf Y(2n)^-$.} The lowest weight 
$\lambda$ is given by:
\be
\begin{array}{l}
    \displaystyle
    u^s\lambda_{i}(u)=\prod_{k=1}^{i}P_{k}(u)\ 
    \prod_{k=i+1}^{n}P_{k}(u+1),\ i=1,\ldots,n\\[2.1ex]
    \displaystyle
    u^s\lambda_{i}(u)=\prod_{k=1}^{n}P_{k}(u+1),\ i=n+1,\ldots,2n
\end{array}
\ee
and the corresponding representation $L(\lambda)$ is built as the 
(subquotient of) the tensor product of 
$s$ evaluation 
representations of $Y(N)$, where $s=dg(\lambda_{i})=\sum_{i=1}^n dg(P_{i})$
and $dg$ is 
the degree in $u$ for $P_{i}$ and in $u^{-1}$ for $\lambda$.

{\bf In the case of $\mathbf Y(N)^+$, $\mathbf N\neq2$:} for $\eps=1$, the construction of $\lambda$
is similar to 
$Y(2n)^-$. For $\eps=3$ ($N=2n$), the construction is still the 
same, but one has to apply the automorphism:
\be
S^{ij}(u)\rightarrow\ S^{i'j'}(u)\ \mbox{ with }\ 
\left\{\begin{array}{l}
i'=N+1-i\ \mbox{ when }\ i=n,n+1\\
i'=i\ \mbox{ otherwise}
\end{array}\right.\label{diese}
\ee
For $\eps=2$, one still considers 
 $L(\lambda)$  the lowest weight representation of $Y(N)$ of lowest weight
$\lambda(u)$ but tensors it with $V_0$ the irreducible representation of $o(N)$ with lowest
weight $(-\half,\ldots,-\half)$. $V_{0}$ can be seen as an evaluation 
representation of the twisted Yangian, where the twisted Yangian is 
embedded into $\cu(o(N))$ using the algebra homomorphism\footnote{Note that the same type 
of homomorphism exists also in the case of $Y^-(2n)$ but is not 
used for the classification of its irreducible finite-dimensional 
representations.}:
\be
\begin{array}{l}
Y^+(N)\ \rightarrow\ \cu(o(N))\\[1.2ex]
\displaystyle 
S(u)\ \rightarrow\ \II+\frac{1}{(u+\half)}\FF
\end{array}
\label{evYtw}
\ee
where $\FF=\sum_{i,j}O^{ij}F_{ij}$, with $O^{ij}$  the 
$o(N)$-generators and $F_{ij}=E_{ij}+E_{N+1-j,N+1-i}$. 

Finally, the case $\eps=4$ is similar 
to $\eps=3$, but with the use of the automorphism (\ref{diese}).

\null

In all the above representations, it may appear that the tensor 
product is not irreducible: in that case, one has to select the 
irreducible part of the lowest weight submodule.

\null

{\bf For $\mathbf Y(2)^+$,} one constructs the lowest weight representation has the 
tensor product of $s$ evaluation representations of $Y(N)$ (where $s$ is the degree of 
$P(u)$), and a 
$o(2)$-representation $V(\eps)$ of weight $\eps-\half$ (understood as 
an evaluation representation of $Y(2)^+$, as in (\ref{evYtw})). 

\null

{\bf Remark:}
Contrarily to the Yangian case, in (\ref{evYtw}), 
the all modes of $S(u)$ have non-zero 
representation. Indeed, while the commutation relation (\ref{com}) is 
invariant under the changes $S(u)\rightarrow\, g(u)S(u)$, the 
symmetry relation (\ref{tauS}) restricts $g(u)$ to be even. Thus, the 
factor $u+\half$ cannot be removed in (\ref{evYtw}).

\sect{Folded \cw-algebras revisited\label{sectWfold}}
It is well-known that the $gl(N)$ Lie algebra can be folded (using an outer
automorphism) into orthogonal and symplectic algebras. In the same way,
folded \cw-algebras have been defined\footnote{Strictly speaking, it is the
folding of "affine" \cw-algebras that has been defined in \cite{fold}, but the
folding of finite \cw-algebras can be defined by the same procedure.}
 in \cite{fold}, and shown to be 
\cw-algebras based on orthogonal and symplectic algebras.

We present here a different proof of this property, adapted to our purpose, 
and generalized to the case of the automorphisms presented in section 
\ref{sectYtw}.
Indeed, we will see that for $N=2n$ one can obtain directly the \cw-algebras
based on $so(2n)$ (as well as those based $sp(2n)$), although there is only
one outer automorphism on $gl(2n)$. The situation is exactly the same
as the one encountered in the folding of $gl(2n)$: although the folding of this
latter algebra (using the automorphism of its Dynkin diagramm) leads to $sp(2n)$, 
it is well-known
that the $so(2n)$ algebra can be constructed as the skew-symmetric matrices of
$gl(2n)$.

\subsection{Automorphism of $gl(Np)$ and $\cw_p(N)$}
As for the Yangian, one introduces an automorphism of $gl(Np)$ defined by
\be
\tau(J_{jm}^{ab})=(-1)^{j+1} \theta^a\theta^b\ J_{jm}^{N+1-b,N+1-a}
\label{tau-glNp}
\ee
where $\theta^a$ is defined as in section \ref{sectYtw}. 
As for the twisted Yangian, there are 
only two different cases to be considered: $\theta_0=1$, or 
$\theta_0=\pm1$ and $N=2n$.

To prove that
$\tau$ is an automorphism of $gl(Np)$, we need the following property
of the Clebsch-Gordan coefficients
\begin{prop}
\label{clebsch}
The Clebsch-Gordan like coefficients obey the rule:
\be
<j,m;p,q|r,s>=(-1)^{j+p+r}\ <p,q;j,m|r,s>
\ee
\end{prop}
\prf
We first prove the property for $j=m$ and $r=s$. 
Due to the property of Clebsch-Gordan coefficient, we have  $q=r-j$ and
we will assume that $q\geq 0$, the proof being similar for $q\leq0$.
One then computes
\beano
<j,j;p,q|r,r> &=&\frac{(-)^r}{\eta_r}\ \sum_{k,\ell,m} 
a^k_{jj}a^l_{pq}a^m_{r,-r}\, tr(E_{k,k+j}E_{l,l+q}E_{m+r,m})\\
&=& \frac{(-)^r}{\eta_r}a^0_{r,-r}\delta_{j+q,r}\sum_{k=1}^{N-j-q} 
a^k_{jj}a^{k+j}_{pq}\\
<p,q;j,j|r,r> &=&\frac{(-)^r}{\eta_r}a^0_{r,-r}\delta_{j+q,r}
\sum_{k=1}^{N-j-q} a^{k+q}_{jj}a^{k}_{pq}\\
&=&\frac{(-)^r}{\eta_r}a^0_{r,-r}\delta_{j+q,r}\sum_{k=1}^{N-j-q} 
(-1)^{2j+p+q} a^{N+1-k-j-q}_{jj}a^{N+1-q-k}_{pq}\\
&=&\frac{(-)^r}{\eta_r}a^0_{r,-r}\delta_{j+q,r}(-1)^{p+q}\sum_{k=1}^{N-j-q} 
 a^{k}_{jj}a^{k+j}_{pq}\\
&=& (-1)^{j+p+r} <j,j;p,q|r,r>
\eeano
which proves the property for $<j,j;p,q|r,r>$.
Now, using
\beano
ad_-(M_{jm}^{ab}) &\equiv& [e_-,M_{jm}^{ab}]=M_{j,m-1}^{ab}\ \Rightarrow\ 
M_{jm}^{ab}=ad_-^{j-m}(M_{jj}^{ab})\\
ad_+(M_{jm}^{ab}) &\equiv& [e_+,M_{jm}^{ab}]=\frac{j(j+1)-m(m+1)}{2}
M_{j,m+1}^{ab}\ \Rightarrow\ 
M_{jm}^{ab}=A_{jm}\ ad_+^{j+m}(M_{jj}^{ab})\\
A_{jm} &=& \prod_{i=-j}^{m-1}\frac{j(j+1)-i(i+1)}{2}=
\frac{(2j)!\, (j+m)!}{2^{j+m}\, (j-m)!}
\eeano
one obtains  the two relations
\beano
<j,m;p,q|r,r> &=& (-1)^{j+m}<j,j;p,q+m-j|r,r>\\
<j,m;p,q|r,s> &=& A_{r,-s}\sum_{i=0}^{r-s} {r-s\choose i}
B^i_{j,m} B^{r-s-i}_{pq}<j,m+i;p,q+r-s-i|r,r>
\eeano
where 
\be
B^i_{j,m} = \prod_{k=1}^{i}\frac{j(j+1)-(m+k)(m+k+1)}{2}=
2^{-i}\frac{(j-m-1)!\,(j+m+i+1)!}{(j-m-i-1)!\,(j+m+1)!}
\ee
These two relations ensure that the property \ref{clebsch} is valid for 
all $<j,m;p,q|r,s>$ coefficients.
\finprf
With this property, it is a simple matter of calculation to show that 
$\tau$ defined in (\ref{tau-glNp}) is an automorphism of $gl(Np)$.
\subsection{Folding $gl(Np)$ and $\cw_p(N)$}
\subsubsection{$gl(Np)$}
One considers the subalgebra Ker$(id-\tau)$ in $gl(Np)$. It is generated
by the combinations:
\be
K_{jm}^{ab}=J_{jm}^{ab}+\tau(J_{jm}^{ab})
\ee
which obey the commutation relations:
\bea
\{K_{jm}^{ab},K_{k\ell}^{cd}\} &=&
\sum_{r=\vert j-k\vert}^{ j+k}\sum_{s=-r}^r\ <j,m;k,\ell\vert r,s>
\,\left(
    \rule{0ex}{2.64ex}
    \delta^{bc}  K^{ad}_{r,s} - 
    \delta^{ad} (-1)^{j+k+r} K^{cb}_{r,s} +\right.\\
&& \left.\rule{0ex}{2.64ex}\hspace{1.2em}
 +(-1)^{j+r}\theta^c\theta^d\delta_{N+1}^{a+c}  K^{N+1-d,b}_{r,s} - 
    \theta^a\theta^b\delta_{N+1}^{b+d} (-1)^{k+r} K^{a,N+1-c}_{r,s}
\right)
\eea
After a rescaling of $K_{jm}^{ab}$ similar to the one given in theorem 
\ref{theo-eps}, one
 recognizes the algebra $so[(2n+1)p]$ (resp. $sp(2np)$; resp. $so(2np)$) when
$N=2n+1$ (resp. $N=2n$, $\theta_0=-1$ ; resp. $N=2n$, $\theta_0=1$).

\null

Looking at the decomposition of the fundamental of $gl(Np)$ with
respect to the principal embedding of $sl(2)$ in $N.sl(p)$ 
(see \cite{classW} and \cite{fold} for the technic used here) one
shows that the subalgebra $N.sl(p)$, generated by the $J^{aa}_{jm}$'s, 
is folded into a $n.sl(p)$ (resp. $n.sl(p)\oplus so(k)$, where $p=2k+1$ is 
chosen odd to get $N$ odd) 
when $N=2np$ (resp. $N=(2n+1)p$).

In the following, we will denote this subalgebra $[N.sl(p)]^\tau$.

\subsubsection{$\cw_p(N)$}
The situation is more delicate, because we are now dealing with the
enveloping algebra of $gl(Np)$, that we denote $\cu[gl(Np)]\equiv\cu(Np)$. 
One introduces the coset:
\[\begin{array}{l}
\cu(Np)^\tau\equiv\cu(Np)/\ck\ \mbox{ where } \ck=\cu(Np)\cdot\cl
\ \mbox{ ; }\cl\mbox{ generated by }
J_{jm}^{ab}-\tau(J_{jm}^{ab}),\ \forall\ a,b,j,m\\
\cw_p(N)^\tau\equiv\cw_p(N)/\cj\ \mbox{ where } \cj=\cw_{p}(N)\cdot\ci 
\ \mbox{ ; }\ci\mbox{ generated by }
W_{j}^{ab}-\tau(W_{j}^{ab}),\ \forall\ a,b,j
\end{array}
\]
We have the property
\begin{prop}
$\tau$ is an automorphism of $\cu(Np)$ 
provided with the Dirac brackets:
\be
\tau\left(\{J_{jm}^{ab},J_{kl}^{cd}\}_*\right)=
\{\tau(J_{jm}^{ab}),\tau(J_{kl}^{cd})\}_*
\ee
Hence, $\tau$ is also an automorphism of $\cw_{p}(N)$.
\end{prop}
\prf
It is obvious that $\tau$ is an automorphism of Poisson brackets on $\cu(Np)$.
Moreover, due to the form of the constraints (\ref{const}), $\tau$ acts as
a relabeling (up to a sign) of the constraints:
\be
\tau(\vph_\alpha)=\vph_{\alpha'}\mbox{ where }\alpha'\equiv\tau(\alpha)
\ee
and also
\be
\tau(\cc_{\alpha\beta})=\cc_{\alpha'\beta'}\ \Rightarrow\ 
\tau(\cc^{\alpha\beta})=\cc^{\alpha'\beta'}
\ee
This shows that this automorphism is compatible with the set 
of constraints $\Phi$ and thus $\tau$ is an automorphism of the Dirac 
brackets. 
\finprf
\begin{coro}
   The Dirac brackets provide $\cw_p(N)^\tau$ with an algebraic structure.
\end{coro}
\prf
We define on $\cw_p(N)^\tau$ a bracket which is just  
the previous Dirac bracket restricted to this coset. Since
$\cw_p(N)^\tau$ is generated by generators of the form $W+\tau(W)$, 
we have:
\beano
\{W+\tau(W),W'+\tau(W')\}_*&=&\{W,W'\}_*+\{\tau(W),\tau(W')\}_*+
\{\tau(W),W'\}_*+\{W,\tau(W')\}_*\\
&=&\{W,W'\}_*+\{\tau(W),W'\}_*+\tau\left(\rule{0ex}{1.2em}\{W,W'\}_*+
\{\tau(W),W'\}_*\right)
\eeano
\finprf
Indeed we have:
\begin{prop}
The $\tau$-folded algebra $\cw_p(N)^\tau$ is the 
$\cw[so[(2n+1)p],n.sl(p)\oplus so(k)]$  algebra
(resp. $\cw[sp(2np),n.sl(p)]$; resp. $\cw[so(2np),n.sl(p)]$ ones)
when $N=2n+1$ and $p=2k+1$ (resp. $N=2n$, $\theta_0=-1$; resp. $N=2n$, 
$\theta_0=1$).
\end{prop}
\prf
On the coset, we have 
$J_{jm}^{ab}\equiv\tau(J_{jm}^{ab})\equiv 2K_{jm}^{ab}$.
We introduce on $\cu(Np)$
\be
2D\vph_\alpha=\vph_\alpha-\vph_{\alpha'}=\vph_\alpha-\tau(\vph_\alpha)
\ ;\ 
2S\vph_\alpha=\vph_\alpha+\vph_{\alpha'}=\vph_\alpha+\tau(\vph_\alpha)
\ee
Since these generators satisfy $D\vph_\alpha=-D\vph_{\alpha'}$ and
$S\vph_\alpha=S\vph_{\alpha'}$ and are in $gl(Np)$, we have
\be
\{S\vph_\alpha,D\vph_\beta\}\in\ci
\ \ie\ \{S\vph_\alpha,D\vph_\beta\}=0\ \mbox{ on }\cw_p(N)^\tau
\ee
Similarly we define
\be
D\cc_{\alpha\beta}=\{D\vph_\alpha,D\vph_\beta\}\ ;\ 
S\cc_{\alpha\beta}=\{S\vph_\alpha,S\vph_\beta\}\ ;\ 
\ee
which obey the properties:
\bea
&&D\cc_{\alpha\beta}=D\cc_{\alpha'\beta'}=-D\cc_{\alpha'\beta}
=-D\cc_{\alpha\beta'} \label{Asym}\\
&&S\cc_{\alpha\beta}=S\cc_{\alpha'\beta'}=S\cc_{\alpha'\beta}
=S\cc_{\alpha\beta'}\label{sym}\\
&&\cc_{\alpha\beta}=S\cc_{\alpha\beta}+
D\cc_{\alpha\beta}\ \mbox{ on }\cw_p^\tau(N)
\eea
We say that a matrix is $\tau$-antisymmetric when it obeys
a relation (\ref{Asym}), and $\tau$-symmetric when it satisfies 
(\ref{sym}). 
$\tau$-antisymmetric matrices are orthogonal to 
$\tau$-symmetric ones:
\[
D\cc\cdot S\cc=0\mbox{ since }
(D\cc\cdot S\cc)_{\alpha\beta}=\sum_\gamma\,
D\cc_{\alpha\gamma}S\cc_{\gamma\beta}=\sum_{\gamma'}\,
D\cc_{\alpha\gamma'}S\cc_{\gamma'\beta}=-\sum_\gamma\,
D\cc_{\alpha\gamma}S\cc_{\gamma\beta}
\]

$S\cc_{\alpha\beta}$ is the matrix of constraints of $gl(Np)^\tau$ 
reduced with respect to $[N.sl(p)]^\tau$. Thus, it is invertible and the associated
Dirac brackets define the algebra $\cw(gl(Np)^\tau,[N.sl(p)]^\tau)$. 
It remains
to show that, on $\cw_p(N)^\tau$, the previously defined Dirac brackets 
coincide with these latter Dirac brackets. 

For that purpose, we use the form $\cc=\cc_0(\II+\wh{\cc})$, given
in \cite{wrtt}, where $\cc_0$ is an invertible $\tau$-symmetric matrix and
$\wh{\cc}$ is nilpotent (of finite order $r$). 
Introducing the $\tau$-symmetrized and antisymmetrized
part of $\wh{\cc}$, one deduces
\be
\cc^{-1}=\cc_0^{-1}\sum_{n=0}^r(-1)^n(S\wh{\cc}+D\wh{\cc})^n
=\cc_0^{-1}\sum_{n=0}^r(-1)^n\left((S\wh{\cc})^n+(D\wh{\cc})^n\right)
=S\cc^{-1}+D\cc^{-1} 
\ee
which shows that $D\cc$ is also invertible. 

On $\cw_p(N)^\tau$, we have
\bea
\{K^{ab}_{(m)},K^{cd}_{(n)}\}_* &=& \{K^{ab}_{(m)},K^{cd}_{(n)}\}-
\{K^{ab}_{(m)},D\vph_\alpha+S\vph_\alpha\}\cc^{\alpha\beta}
\{D\vph_\beta+S\vph_\beta,K^{cd}_{(n)}\}\\
&=&\{K^{ab}_{(m)},K^{cd}_{(n)}\}-
\{K^{ab}_{(m)},S\vph_\alpha\}\cc^{\alpha\beta}
\{S\vph_\beta,K^{cd}_{(n)}\}\\
&=& \{K^{ab}_{(m)},K^{cd}_{(n)}\}-
\{K^{ab}_{(m)},S\vph_\alpha\}(S\cc^{\alpha\beta}+D\cc^{\alpha\beta})
\{S\vph_\beta,K^{cd}_{(n)}\}
\eea
From the $\tau$-antisymmetry of $D\cc^{-1}$, we get
\be
\{.,S\vph_\alpha\}D\cc^{\alpha\beta}\{S\vph_\beta,.\}=
\{.,S\vph_{\alpha'}\}D\cc^{\alpha'\beta}\{S\vph_\beta,.\}=
-\{.,S\vph_\alpha\}D\cc^{\alpha\beta}\{S\vph_\beta,.\}=0
\ee
which leads to the Dirac brackets:
\be
\{K^{ab}_{(m)},K^{cd}_{(n)}\}_* = \{K^{ab}_{(m)},K^{cd}_{(n)}\}-
\{K^{ab}_{(m)},S\vph_\alpha\}S\cc^{\alpha\beta}
\{S\vph_\beta,K^{cd}_{(n)}\}
\ee
These Dirac brackets are the $\cw(gl(Np)^\tau,[N.sl(p)]^\tau)$
algebra ones, by definition of $S\cc$.
\finprf

\sect{Folded \cw-algebras as truncated twisted Yangians\label{results}}
We consider the $\cw_p(N)$ algebra in the Yangian basis. The Poisson
brackets are
\be
\{T_{(q)1},T_{(r)2}\} = \sum_{s=0}^{min(p,q,r)-1} (
P_{12}T_{(s)1}T_{(q+r-s)2}-T_{(s)2}T_{(q+r-s)1}P_{12})
\label{zozo}
\ee
with the convention $T_{(r)}=0$ for $r>p$. The action of the automorphism
$\tau$, both for twisted Yangian and folded $\cw_p(N)$ algebra, reads
\be
\tau(T_{(m)})=(-1)^{m}T^t_{(m)}
\ee
However, from the twisted Yangian point of view, one selects
 the generators 
 \[
 S_{(m)}=\sum_{r+s=m}(-1)^{s}T_{(r)}T^t_{(s)}
 \]
while in the folded \cw-algebra case, one constrains the generators to
$T_{(m)}=(-1)^{m}T^t_{(m)}$. Although the procedures are different 
(and indeed lead to different generators), we have:
\begin{theor}
As an algebra, the folded \cw-algebra $\cw(gl(Np)^\tau,[N.sl(p)]^\tau)$ is
isomorphic to the truncation (at level $p$) of the 
(classical) twisted Yangian $Y(N)^\tau$. 

More precisely, we have the correspondences:
\be
\begin{array}{ccl}
Y_{p}(2n)^- &\longleftrightarrow & \cw[sp(2np), n.sl(p)] \\
Y_{p}(2n)^+ &\longleftrightarrow & \cw[so(2np), n.sl(p)] \\
Y_{p}(2n+1)^+ &\longleftrightarrow & \cw[so((2n+1)p), n.sl(p)\oplus 
so(k)] \ ;\ p=2k+1\\
\end{array}
\ee
\end{theor}
\prf
We prove this theorem by showing that the Dirac brackets of the folded 
\cw-algebra coincide with the Poisson brackets 
\be
\{S_{1}(u),S_{2}(v)\} = [r_{12}(u-v), S_{1}(u)S_{2}(v)] +
S_{2}(v)r'_{12}(u+v)S_{1}(u)-S_{1}(u)r'_{12}(u+v)S_{2}(v)
\label{truc}
\ee
with the truncation $S_{(m)}=0$ for $m>p$. 

We start with the $\cw_p(N)$ algebra in the truncated Yangian basis and define 
\be
2\vph_{(s)}= T_{(s)}-(-1)^s\, T^t_{(s)} \ \mbox{ and }\ 
2K_{(s)}=T_{(s)}+(-1)^s\, T^t_{(s)}
\ee
The folding (of the \cw-algebra) corresponds to
\be
\vph_{(s)}=0\ \ie K_{(s)}=(-1)^sK^t_{(s)}\label{eq:fold}
\ee
It is a simple matter of calculation using (\ref{zozo}), to compute
\bea
2\{K_{(q)1},K_{(r)2}\} &=& \sum_{s=0}^{M}\left[ \rule{0ex}{2.4ex}
P_{12}K_{(s)1}K_{(r+q-s-1)2}-
K_{(r+q-s-1)2}K_{(s)1}P_{12}+\right.\nonu
&&\left.+(-1)^{q+s}\left(\rule{0ex}{2.1ex}K_{(s)1}Q_{12}K_{(r+q-s-1)2}-
K_{(r+q-s-1)2}Q_{12}K_{(s)1}\right)\right]
\eea
which is equivalent to the relation (\ref{truc}) for $S(u)\equiv K(\frac{u}{2})$. 
The constraint (\ref{eq:fold}) is then rewritten as
$S^t(-u)=S(u)$. Thus, the folded \cw-algebra and the truncated 
twisted Yangian are defined by the same relations.
\finprf

\subsection{Quantization of \cw-algebras}
Now that folded \cw-algebras have proved to be truncation of twisted 
Yangians, there quantization is very simple. It takes the form
\bea
&&R_{12}(u-v)\, S_{1}(u)\, R'_{12}(u+v)\, S_{2}(v) = 
S_{2}(v)\, R'_{12}(u+v)\, S_{1}(u)\, R_{12}(u-v)\\[2.1ex]
&&\mbox{with }\left\{\begin{array}{l}\displaystyle
 R_{12}(x)=\II\otimes\II-\frac{1}{x}P_{12}\ ;\ R'_{12}(x)=\II\otimes\II
 -\frac{1}{x}Q_{12}\\
\displaystyle S(u)=\sum_{m=0}^{p}\, u^{-m}\, S_{(m)}\ ;\ S_{(0)}=\II
\end{array}\right.
\eea
Let us remark that, contrarily to the quantization of $\cw_{p}(N)$, 
the quantization of $\cw_{p}(N)^\tau$ indeed modify the commutation 
relations, adding a non-trivial non-central new term (see equation 
(\ref{PB-S}) with respect to (\ref{com-S}), its quantization).

\subsection{Center and finite-dimensional irreducible representations}
Starting from the classification of finite dimensional representations of the twisted 
Yangian, one can deduce the ones of the folded \cw-algebras. 
\begin{prop}
The finite dimensional irreducible representations  of 
$\cw_p(N)^\tau\neq \cw_{p}^+(2)$ are 
given by the property \ref{irrepsYtw} with the restriction that the polynomials
$P_i(u)$  and the parameter $\eps$ must obey to the following constraints
\be
\sum_{i=1}^{\bn} dg(P_i)\leq \frac{p}{2} \ \mbox{ and } \eps\mbox{ odd }
\ee
where $dg(P_i)$ is the degree of $P_i(u)$. 

In the special case $\cw_{p}^+(2)\equiv\cw[so(2p),sl(p)]$, 
one imposes $dg(P)\leq \frac{p}{2}$ and $\eps=\half$.
\end{prop}
\prf
The calculation is the same as for $\cw_p(N)$, starting from $Y(N)$ 
(see \cite{wrtt} for more detail). It essentially relies on the fact that 
in the tensor product $k$ evaluation representations (of the Yangian 
$Y(N)$), we have 
\[
\pi_{k}(T_{(m)})=0\ \Leftrightarrow\ m>k
\]
where $\pi_{k}$ is the representation morphism. This property is also 
valid when considering the subquotient of the tensor product. 
 In the case of $Y(N)$, this number corresponds to the sum of the 
 degrees of the polynomials $P_{i}$, hence the condition $k\leq p$ to get a representation of 
$Y_{p}(N)$.
In the case of $Y^\pm(N)$, considering tensor products of $Y(N)$-evaluation 
representations, and since $S(u)$ is quadratic in $T(u)$, one gets:
\[
\pi_{k}(S_{(m)})=0\ \Leftrightarrow\ m>2k
\]
Hence, the sum of the degree of the polynomials 
$P_{i}$ is (up to the $V_{0}$ representation)
half of the number of $Y(N)$-evaluation representations 
used to build the $Y^\pm(N)$ representation (see constructions in 
proofs of theorem 5.8, and followings in \cite{Mol}). When an 
$o(N)$-representation $V_{0}$ (or a $o(2)$-representation $V(\eps)$ for $Y^+(2)$)
is involved, all the generators have non-vanishing representation, 
and thus cannot be set to zero (see remark at the end of section 
\ref{sectYtw}): we have not a representation of the 
truncated Yangian. Hence, only the values $\eps=1,3$ are allowed for 
$Y^\pm(N)\neq Y^+(2)$, and only the value $\eps=\half$ (\ie $V(\eps)$ 
trivial representation) for $Y^+(2)$.

Conversely, starting with a finite dimensional irreducible 
representation of the truncated twisted Yangian, one can construct a 
representation of the whole twisted Yangian by representing the 
remaining generators by zero. This representation is obviously 
irreducible (since it is for the truncated twisted Yangian), and thus 
falls into the classification of theorem \ref{irrepsYtw}.
\finprf

Again, we can follow  the same steps as in \cite{wrtt} to conclude:
\begin{prop}
The center of $\cw_p(N)^\tau$ has dimension $\bn\, p$ and is generated by the
$\bn\, p$ first even Casimir operators of the underlying $gl(Np)$ algebra.
\end{prop}
\prf
From the Hamiltonian reduction on $[gl(Np)]^\tau$, we known that the 
center of $\cw_{p}(N)^\tau$ contains the Casimir operators of 
$[gl(Np)]^\tau$ and that the $\bn\, p$ first ones are algebraically 
independent. The property (\ref{center}) and the truncation show that 
they are the only ones.
\finprf

\sect{Conclusion\label{concl}}
We have shown that the truncation of twisted Yangians $Y^\pm(N)$ are isomorphic
to finite \cw-algebras based on orthogonal or symplectic algebras. As for $Y(N)$ and
finite \cw-algebras based on $gl(N)$, this isomorphism allows to classify
all the finite dimensional irreducible representations of these 
\cw-algebras, and to determine their center. 
It provides also a $R$-matrix formulation of the \cw-algebras. However, contrarily
to the case of $gl(N)$, the formulation is an $ABCD$-type one. This confirms the
remark already done for $\cw_p(N)$ that there seems to be no natural Hopf 
structure on \cw-algebras. 

On the other hand, the fact that the isomorphism between Yangians and 
finite \cw-algebras is still valid for $so(m)$ and $sp(m)$ algebras is 
a good point in favor of the generalization of Yangians. Indeed works 
are in progress for an $R$-matrix for all the finite \cw-algebras, and 
a limiting procedure on it should lead to such  generalized Yangians.

{F}inally, note that the supersymmetrization of these 
construction can also be done \cite{superW}.

\null

{\bf Acknowledgments}

\null

We warmly thank D. Arnaudon for his patience, and all the fruitful 
implied discussions.

\end{document}